\documentclass[12pt]{article}
\usepackage{amsmath}
\usepackage{amssymb}
\newtheorem{thm}{Theorem}
\newtheorem{lemma}[thm]{Lemma}
\newtheorem{cor}[thm]{Corollary}

\newtheorem{definition}[thm]{Definition}
\DeclareMathOperator{\divergence}{div}

\begin{document}
\title{On the Time Derivative in an Obstacle Problem}
\author{Peter Lindvist}
\date{Norwegian University of Science and Technology}
\maketitle

\maketitle

{\small \textsc{Abstract:} \textsf{We prove that the time derivative
 of the solution for
    the obstacle problem related to the Evolutionary $p$-Laplace Equation
    exists in Sobolev's sense, provided that the given obstacle is
 smooth enough. We keep $p \geq 2.$}}
\section{Introduction}

The celebrated \emph{Evolutionary $p$-Laplace Equation} is much
studied and the regularity theory for the solutions is almost
complete.
 We refer to the
book [dB] about this fascinating equation.  In general, the
corresponding
 subsolutions and supersolutions do not
possess that much regularity, they are semicontinuous. 
We are interested in a special kind of weak supersolutions of the 
Evolutionary $p$-Laplace Equation, namely the solutions of an
obstacle problem. In the presence of a smooth obstacle the regularity
improves a lot. Given a function $\psi = \psi(x,t)$ in a bounded domain
$\Omega_T = \Omega \times (0,T),$ where $\Omega \subset \mathbf{R}^n,$
 we consider all functions $v$ such that
$$ \frac{\partial v}{\partial t} \geq \nabla\! \cdot \! (|\nabla
v|^{p-2}\nabla v) \quad \text{and} \quad v \geq \psi\quad \text{in}
\quad \Omega_T.  $$
 The function $\psi$ acts as an \emph{obstacle}. The
smallest admissible $v$ is the solution of the obstacle
problem. (This makes sense because a comparison principle is valid.)
 However, the above description was only formal. We will
instead use Definition 1 below, which is more adequate since it
comes with a \textsf{variational inequality.} ---We will restrict ourselves to
the case $p > 2$, the so-called slow diffusion case.

It is an established fact that if the obstacle $\psi$ is smooth
enough, the solution to the obstacle problem inherits some
regularity. Our objective is the time derivative $u_t$ of the solution
$u,$ which \emph{a priori} is only known to be a distribution. Our
main result Theorem 2 states that, if $\psi$ has continuous second
derivatives, then the time derivative $u_t$ exists in Sobolev's sense
and  it belongs to the space $L^{p/(p-1)}_{loc}(\Omega_T).$ A formula is
given for the derivative. The most laborious part of the proof is to
show that $\Delta_p u = \nabla \cdot (|\nabla v|^{p-2}\nabla v)$ is a function
so that the rule
\begin{equation*}
\int_0^T\!\!\int_{\Omega}\varphi\Delta_p u \,dx\,dt =
  -\int_0^T\!\!\int_{\Omega}\langle \, |\nabla
  u|^{p-2}\nabla u,\nabla \varphi \rangle \,dx\,dt\
\end{equation*}
with test functions applies. The equation has first to be
regularized, keeping the obstacle unaffected, and then 
difference quotients are used. The test functions in [L1] can be
adjusted to work here.

An important feature, typical for obstacle problems, is that in the open
set $\Upsilon =\{u > \psi\}$ where the obstacle does not hinder, $u$ is,
actually, a solution to the differential equation. Thus in $\Upsilon$
the equation  $u_t = \Delta_pu$ holds in the weak sense. The boundary
of the \emph{coincidence set} $\Xi = \{u=\psi\}$ is crucial. This
enables us to get an identity for the integral $\iint
u\varphi_t\,dx\,dt,$
from which one can deduce the existence of the time derivative sought
for. The special case with no obstacle present was treated in [L2]. ---See
also [BDM] for some general comments valid for ``irregular'' obstacles.

{\small To this we may add a curious fact valid for $\psi \in
C^{2}(\Omega_T).$ At all points in the coincidence
set $\Xi$ the obstacle satisfies the inequality
$$\frac{\partial \psi}{\partial t} \geq \Delta_p\psi.$$
 Thus a point at which $\frac{\partial \psi}{\partial t} <
\Delta_p\psi$ cannot belong to the coincidence set. This piece of
information follows from the characterization of continuous
supersolutions as \emph{viscosity} supersolutions, cf. [JLM].
  Then $\psi$ itself
can do as a test function for the pointwise testing required in the
theory of viscosity solutions. (The reader may consult [K] for
some basic concepts.) ---We will not need this observation.}

It is likely that the time derivative belongs to the space
$L^{2}_{loc}(\Omega_T)$, but an eventual proof of this improvement
would require much stronger regularity considerations for $\nabla u.$
 We have kept $p > 2$, but one can  expect a counterpart to Theorem
 2 valid in the  extended range $p > 2n/(n+2).$ The difficulty about
 further generalizations with $\Delta_p u$ replaced by some operator
 $\divergence \mathbf{A}_p$  is the following. It is absolutely
 necessary that the solutions of the differential equation
$$\frac{\partial u}{\partial t} = \divergence
\mathbf{A}_p(x,t,u,\nabla u)$$
enjoy the property of having a time derivative themselves, in order
that the corresponding results could be extended to the related
 obstacle problem. This considerably restricts the possibilities.
\section{Preliminaries}

Let  $\Omega$ be a bounded domain in the $n$-dimensional space $\mathbf{R}^n$
having a Lipschitz regular boundary. Suppose that a function $\psi =
\psi(x,t)$ is given in the closure of the space-time cylinder
$\Omega_T = \Omega \times (0,T)$. The function $\psi$ acts as an
\emph{obstacle} so that the admissible functions are forced to lie
above $\psi$ in $\Omega_T$. We make the
 \begin{equation*}
\text{\textsf{Assumption:}}\quad\psi \in C(\overline{\Omega_T})
 \cap W^{2,p}(\Omega_T).\end{equation*}
For simplicity the obstacle $\psi$ also determines the values of the
admissible functions on the \emph{parabolic boundary} 
$$\Gamma_T = \Omega \times \{0\}\; \cup  \; \partial \Omega
\times [0,T].$$ The \emph{class of admissible functions} is
\begin{gather*}
\mathcal{F}_{\psi} = \{v\in L^{p}(0,T;W^{1,p}(\Omega))|\,v \in
C(\overline{\Omega_T}),\: v \geq \psi \ \text{in}\ 
 \Omega_T,\: v = \psi\ \text{on}\ \Gamma_T \}.
\end{gather*}
---We keep $p \geq 2.$

\begin{definition} We say that the function $u \in \mathcal{F}_{\psi}$ is the
  solution to the obstacle problem, if the inequality
\begin{gather}
\int_0^T\!\!\int_{\Omega}\Bigl(\langle|\nabla u|^{p-2}\nabla u,\nabla(\phi
-u)\rangle + (\phi - u)\,\frac{\partial \phi}{\partial t}\Bigr)dx\,dt
\nonumber \\
\geq \frac{1}{2}\int_{\Omega}(\phi(x,T)-u(x,T))^{2}\,dx
\end{gather}
holds for all \emph{smooth} functions $\phi \in  \mathcal{F}_{\psi}.$
\end{definition}
The solution exists and is unique, cf.[AL] and [C]. See also [KKS].
 It is also a supersolution of the
equation $u_t \geq \Delta_p u$, i.e.,
\begin{equation}
\int_0^T\!\!\int_{\Omega}\Bigl(\langle|\nabla u|^{p-2}\nabla u,\nabla
\varphi \rangle - u\,\frac{\partial \varphi}{\partial t}\Bigr)dx\,dt
\geq 0
\end{equation}
for all non-negative $\varphi \in C_0^{\infty}(\Omega_T).$ 

Notice that nothing is assumed about the time derivative $u_t$. Our
main result is the theorem below.
\begin{thm} The time derivative $u_t$ of the solution $u$ to the
  obstacle problem exists in the Sobolev sense and $u_t \in
  L^{p/(p-1)}_{loc}(\Omega_T).$ It is the function 
\begin{equation*}
u_t = \begin{cases}
\psi_t\quad \text{in}\quad \Xi\\
\Delta_p u \quad\text{in}\quad \Omega_T \setminus \Xi
\end{cases}
\end{equation*}
where $\Xi = \{u = \psi\}$ denotes the coincidence set.
\end{thm}

 In
order to avoid the difficulty with the ``forbidden'' time derivative
$u_t$ in the proof, we have to regularize the equation, keeping the obstacle
unchanged. We replace $|\nabla u|^{p-2}\nabla u$ by 
$$\Bigl(|\nabla u|^{2} + \varepsilon^{2}\Bigr)^{\frac{p-2}{2}}\nabla u$$
to obtain an equation which does not degenerate as $\nabla u = 0.$

\begin{lemma} There is a unique $u^{\varepsilon} \in
  \mathcal{F}_{\psi}$ such that
\begin{gather}
\int_0^T\!\!\int_{\Omega}\Bigl(\langle\bigl|\nabla u^{\varepsilon}|^{2}
+ \varepsilon^{2}\bigr)^{\frac{p-2}{2}}\nabla u^{\varepsilon},\nabla(\phi
-u^{\varepsilon})\rangle + (\phi - u^{\varepsilon})\,\frac{\partial \phi}{\partial t}\Bigr)dx\,dt \nonumber\\
\geq \frac{1}{2}\int_{\Omega}(\phi(x,T)-u^{\varepsilon}(x,T))^{2}\,dx
\end{gather}
for all smooth functions $\phi$ in the class $\mathcal{F}_{\psi}.$ In
the open set $\{u^{\varepsilon} > \psi \}$ the function
$u^{\varepsilon}$ is a solution of the equation $$u^{\varepsilon}_t =
\nabla\! \cdot\! \Bigl(\bigl( |\nabla u^{\varepsilon}|^{2}
+ \varepsilon^{2}\bigr)^{\frac{p-2}{2}}\nabla
u^{\varepsilon}\Bigr).$$ In the case $\varepsilon \not = 0$ we have
$u^{\varepsilon} \in C^{\infty}(\Omega_T)$ and $\frac{\partial
  u^{\varepsilon}}{\partial t} \in L^{2}(\Omega_T).$
\end{lemma}

\emph{Proof:} The existence can be extracted from the proof of [AL,
Theorem 3.2]. The regularity for the nondegenerate case $\varepsilon
\not = 0$ is according to the standard parabolic theory described in
the celebrated book [LSU]. The proof of the H\"{o}lder continuity for
the degenerate case $\varepsilon = 0$ is in [C]. 

\medskip
When $\varepsilon \not = 0$, we can rewrite equation (3) in the more
convenient form
\begin{gather}
\int_0^T\!\!\int_{\Omega}\Bigl(\langle\underbrace{\bigl|\nabla
  u^{\varepsilon}|^{2} +
\varepsilon^{2}\bigr)^{\frac{p-2}{2}} 
 \nabla u^{\varepsilon}} _{\mathbf{A}_{\varepsilon}(x,t)}  ,\nabla \eta\rangle
 + \eta\,\frac{\partial
  u^{\varepsilon}  }{\partial t}\Bigr)dx\,dt \geq 0
\end{gather}
valid for all test functions $\eta$ such that $\eta \geq \psi \
 -u^{\varepsilon}$ in $\Omega_T$ and $\eta = 0$ on $\Gamma_T.$ We may
 even use any
 continuous $\eta \in L^{p}(0,T;W_{0}^{1,p}(\Omega))$ with $\eta(x,0) = 0.$

In order to proceed to the limit under the integral sign in the
forthcoming equations we need the convergence result below, where $u$
denotes the solution to the original obstacle problem, the one with
$\varepsilon = 0.$

\begin{lemma}
\begin{equation}
\lim_{k\rightarrow 0}\, \int_0^T\!\!\int_{\Omega}\Bigl(|u^{\varepsilon}
  - u|^p + |\nabla u^{\varepsilon} - \nabla u|^p\Bigr)\,dx\,dt\, =
  \, 0.
\end{equation}
\end{lemma}

\medskip
\emph{Proof:} It was established in [KL, Lemma 3.2] that 
\begin{equation}
\lim_{k\rightarrow 0}\, \int_0^T\!\!\int_{\Omega} |\nabla
u^{\varepsilon} - \nabla u|^p\,dx\,dt\, = \, 0,\end{equation}
but the  strong convergence of the functions themselves requires, as
it were, an extra
compactness argument. Since $u^{\varepsilon}$ is a weak supersolution,
there exists a Radon measure $\mu_{\varepsilon}$ such that 
\begin{gather*}
\int_0^T\!\!\int_{\Omega}\Bigl(\langle\bigl|\nabla u^{\varepsilon}|^{2}
+ \varepsilon^{2}\bigr)^{\frac{p-2}{2}}\nabla u^{\varepsilon},\nabla \varphi
\rangle -  u^{\varepsilon}\,\frac{\partial \varphi}{\partial t}\Bigr)dx\,dt 
 = \int_{\Omega_T}\varphi \,d\mu_{\varepsilon}
\end{gather*}
for all functions $\varphi \in  C^{\infty}_{0}(\Omega_T),$ whether positive or
not. This is a consequence of Riesz's Representation Theorem, cf. [EG,
§1.8]. See
[KLP] for details.

Given a regular open set (for example a polyhedron) $U \subset \subset
\Omega_T,$ we have to verify that
$$\mu_{\varepsilon}(U) \leq M_U$$
with a bound independent of $\varepsilon,\, 0< \varepsilon < 1.$
Then the lemma follows as in [KLP, pp. 720-721]. (There [S] was used.) To
this end, choose a test function $\varphi  \in
C^{\infty}_{0}(\Omega_T)$ such that $0 \leq \varphi \leq 1$ and
$\varphi = 1$ in $U$. A rough estimation yields
\begin{gather*}
\mu_{\varepsilon}(U) = \int_{\Omega_T} \,d\mu_{\varepsilon} =
\int_{\Omega_T}\varphi \,d\mu_{\varepsilon}\\
= \int_0^T\!\!\int_{\Omega}\Bigl(\langle\bigl|\nabla u^{\varepsilon}|^{2}
+ \varepsilon^{2}\bigr)^{\frac{p-2}{2}}\nabla u^{\varepsilon},\nabla \varphi
\rangle -  u^{\varepsilon}\,\frac{\partial \varphi}{\partial
  t}\Bigr)dx\,dt\\
\leq  
C_1\left(\|\nabla u^{\varepsilon}\|^{p}_{L^{p}(\Omega_T)} +
  \varepsilon^{\frac{p(p-2)}{p-1}}\right) + C_2\|u^{\varepsilon}\|_{\infty}.
\end{gather*}
By the maximum principle $\|u^{\varepsilon}\|_{\infty} \leq \|\psi\|_{\infty}$
and $ \|\nabla u^{\varepsilon}\|^{p}_{L^{p}(\Omega_T)}$ is uniformly
bounded, since the gradients converge strongly. This yields the bound $M_U.$
$\Box$
\section{The  gradient estimate}

In order to prove that $\Delta_p u$ is a \emph{function}, $u$ denoting
the solution to the obstacle problem, we show that the function
$\textbf{F}  = |\nabla u|^{(p-2)/2}\nabla u$, where the usual power
$p-2$ has been replaced by $(p-2)/2$, is in a suitable first
order Sobolev $x$-space. This will immediately imply the desired
result. At a first reading one had better to assume that the obstacle
$\psi$ is as smooth as one pleases, say of class
$C^{2}(\overline{\Omega_T}).$
Actually, only the Sobolev derivatives $\psi_{x_i x_j}$ and $\psi_{x_i
  t}$
are needed, while $\psi_{tt}$ does not appear at all. We recall our
assumption
$\psi \in C(\overline{\Omega_T}) \cap W^{2,p}(\Omega_T)$
and use the abbreviation
$$|D^2\psi|^2 = \sum \psi_{x_ix_j}^{2}.$$
Under these assumptions about the obstacle $\psi = \psi(x,t)$ we have the
following result.

\begin{thm}
For the solution $u$ to the obstacle problem,
the derivative $D\textbf{F}$ of $$\textbf{F} = |\nabla
u|^{\frac{(p-2)}{2}} \nabla u$$ exists in Sobolev's sense and belongs to
  $L^{p/(p-1)}_{loc}(\Omega_T)$. The estimate
\begin{gather*}
\int_0^T\!\!\int_{\Omega} \zeta^p|D\textbf{F}|^2\,dx\,dt \leq
C \int_0^T\!\!\int_{\Omega}( \zeta^p + | \nabla \zeta|^p)|\nabla
u|^p\,dx\,dt\\
+ C \int_0^T\!\!\int_{\Omega} \zeta^p |\nabla u|^2\,dx\,dt +
 C \int_0^T\!\!\int_{\Omega}|\nabla \zeta|^p |\nabla
 \psi|^p\,dx\,dt\\
+ C \int_0^T\!\!\int_{\Omega} \zeta^p (|D^2 \psi|^p + |\nabla
 \psi_{t}|^2)\,dx\,dt + C \int_{\Omega} \zeta^p  | \nabla
 \psi(x,T)|^2\,dx
\end{gather*}
holds for each non-negative test function $\zeta = \zeta(x)$ in
$C_0^{\infty}(\Omega)$; and $C = C(p).$
\end{thm}

\emph{Proof:} The proof is based on the  regularized obstacle problem
and equation
(4), where we abbreviate
$$ \mathbf{A}_{\varepsilon}(x,t) = \bigl(|\nabla u^{\varepsilon}|^{2}
+
\varepsilon^{2}\bigr)^{\frac{p-2}{2}} 
 \nabla u^{\varepsilon}.$$
 We denote its solution by $u$, suppressing the index
$\varepsilon$. Thus $u$ means $u^{\varepsilon}$, to begin with. Given
$\zeta$, the variable $x$ is given a small increment $h$ so that the
test function
\begin{gather*}
\eta = \psi(x,t)-u(x,t) + \zeta(x)^p[u(x+h,t)-\psi(x+h,t)]\\
= \zeta(x)^p\overbrace{[u(x+h,t)-u(x,t)]}^{\mathbf{\Delta_h}u}\, - \,
\zeta(x)^p\overbrace{[\psi(x+h,t)-\psi(x,t)]}^{\mathbf{\Delta_h}\psi}
\\
-(1-\zeta(x)^p)[u(x,t)-\psi(x,t)]
\end{gather*}
is admissible in the regularized equation 
\begin{gather}
\int_0^T\!\!\int_{\Omega}\Bigl(\langle\mathbf{A}_{\varepsilon}(x,t),
\nabla \eta\rangle
 + \eta\,\frac{\partial
  u }{\partial t}\Bigr)dx\,dt \geq 0.
\end{gather}
Inserting the test function, we obtain
\begin{align*}
\phantom{+}&\int_0^T\!\!\int_{\Omega}\Bigl(\langle\mathbf{A_{\varepsilon}}(x,t),\nabla
(\zeta^p\mathbf{\Delta_h}u)\rangle 
+\zeta^p\mathbf{\Delta_h}u\frac{\partial u}{\partial
  t}\Bigr)\,dx\,dt\\
-&\int_0^T\!\!\int_{\Omega}\Bigl(\langle\mathbf{A_{\varepsilon}}(x,t),\nabla
(\zeta^p\mathbf{\Delta_h}\psi)\rangle 
+\zeta^p\mathbf{\Delta_h}\psi\frac{\partial u}{\partial
  t}\Bigr)\,dx\,dt\\
\geq &\int_0^T\!\!\int_{\Omega}\Bigl(\langle\mathbf{A_{\varepsilon}}(x,t),\nabla
\bigl((1-\zeta(x)^p)[u(x,t)-\psi(x,t)]\bigr)\rangle \\
&\phantom{abbbbb} +(1-\zeta(x)^p)[u(x,t)-\psi(x,t)]\frac{\partial u}{\partial
  t}\Bigr)\,dx\,dt\\
&\geq 0.
\end{align*}
The last integral is non-negative, because 
\begin{equation*}
(1-\zeta(x)^p)[u(x,t)-\psi(x,t)]
\end{equation*}
will do as a test function in the equation (7). This observation is
important here.

Aiming at difference quotients we give $x$ the increment $h$.
 The translated function $u(x+h,t)$ solves the obstacle problem with
the translated obstacle $\psi(x+h,t)$, all this with respect to the shifted
domain
$\Omega^{h}\times(0,T)$ where $\Omega^{h} = \{x|\,x+h \in \Omega\}$. For
sufficiently small $h$ we have
\begin{equation}
\int_0^T\!\!\int_{\Omega^{h}}\Bigl(\langle\mathbf{A_{\varepsilon}}(x+h,t),\nabla
 \eta(x,t)\rangle 
+\eta(x,t)\,\frac{\partial u(x+h,t)}{\partial
  t}\Bigr)\,dx\,dt \geq 0
\end{equation}
whenever $\eta(x,t) \geq \psi(x+h,t) - u(x+h,t)$ and $\eta = 0$ on the
parabolic boundary of $\Omega^{h}\times(0,T)$. Here
\begin{gather*}
\eta = \psi(x+h,t)-u(x+h,t) + \zeta(x)^p[u(x,t)-\psi(x,t)]\\
= \zeta(x)^p\overbrace{[u(x+h,t)-u(x,t)]}^{\mathbf{\Delta_h}u}\, - \,
\zeta^p(x)\overbrace{[\psi(x+h,t)-\psi(x,t)]}^{\mathbf{\Delta_h}\psi}
\\
-(1-\zeta(x)^p)[u(x+h,t)-\psi(x+h,t)]
\end{gather*}
will do. We obtain
\begin{align*}
-&\int_0^T\!\!\int_{\Omega^{h}}\Bigl(\langle\mathbf{A_{\varepsilon}}(x+h,t),\nabla
(\zeta^p\mathbf{\Delta_{h}}u)\rangle 
+\zeta^p\mathbf{\Delta_h}u\,\frac{\partial u(x+h,t)}{\partial
  t}\Bigr)\,dx\,dt\\
+&\int_0^T\!\!\int_{\Omega^{h}}\Bigl(\langle\mathbf{A_{\varepsilon}}(x+h,t),\nabla
(\zeta^p\mathbf{\Delta_h}\psi)\rangle 
+\zeta^p\mathbf{\Delta_h}\psi\,\frac{\partial u(x+h,t)}{\partial
  t}\Bigr)\,dx\,dt\\
\geq
&\int_0^T\!\!\int_{\Omega^{h}}\Bigl(\langle\mathbf{A_{\varepsilon}}(x+h,t),\nabla 
\bigl((1-\zeta(x)^p)[u(x+h,t)-\psi(x+h,t)]\bigr)\rangle \\
&\phantom{abbbbb}+(1-\zeta(x)^p)[u(x+h,t)-\psi(x+h,t)]\,\frac{\partial
  u(x+h,t)}{\partial t}\Bigr)\,dx\,dt\\
&\geq 0.
\end{align*}
The last integral is positive
because $$(1-\zeta(x)^p)[u(x+h,t)-\psi(x+h,t)]$$ will do as a test
function in the translated equation (8). This observation is
essential here. The integrals in the left-hand member of the
inequality are, in fact, taken only over the support of the
function $\zeta(x)$. Hence we have an inequality with
integrals taken only over $\Omega_T$, provided that $|h| <
\mathrm{dist(supp}\zeta,\partial \Omega)$. Thus $\Omega^{h}$ is no
longer directly involved.

We add the two estimates, grouping the differences, and obtain
\begin{align*}
+&\int_0^T\!\!\int_{\Omega}\langle\mathbf{A_{\varepsilon}}(x+h,t)
-\mathbf{A_{\varepsilon}}(x,t),\nabla
(\zeta^p\mathbf{\Delta_h}u)\rangle 
\,dx\,dt\\
\leq &\int_0^T\!\!\int_{\Omega}\langle\mathbf{A_{\varepsilon}}(x+h,t)
-\mathbf{A_{\varepsilon}}(x,t) ,\nabla
(\zeta^p\mathbf{\Delta_h}\psi)\rangle\,dx\,dt\\
-&\int_0^T\!\!\int_{\Omega}\zeta^p\mathbf{\Delta_h}u \cdot
\mathbf{\Delta_h}\Bigl( \frac{\partial u}{\partial
  t}\Bigr)\,dx\,dt +
\int_0^T\!\!\int_{\Omega}\zeta^p\mathbf{\Delta_h}\psi \cdot 
\mathbf{\Delta_h}\Bigl( \frac{\partial u}{\partial
  t}\Bigr)\,dx\,dt.
\end{align*}
The integrals with the time derivatives can be integrated by parts:
\begin{align*}
-&\int_0^T\!\!\int_{\Omega}\zeta^p  \frac{\partial }{\partial t}
\frac{  \mathbf{(\Delta_h}u)^2 }{2}\,dx\,dt+
\int_0^T\!\!\int_{\Omega}\zeta^p \mathbf{\Delta_h}\psi \cdot 
\mathbf{\Delta_h}\Bigl( \frac{\partial u}{\partial
  t}\Bigr)\,dx\,dt\\
=&-\int_{\Omega}\zeta^p(x)\frac{\mathbf{(\Delta_h}u)^2 }{2}\Big\lvert_0^T\,dx +
\int_{\Omega}\zeta^p(x)\mathbf{\Delta_h}\psi\cdot
\mathbf{\Delta_h}u\Big\rvert_0^T\,dx \\
-&\int_0^T\!\!\int_{\Omega}\zeta^p\mathbf{\Delta_h}u \cdot 
\mathbf{\Delta_h}\Bigl( \frac{\partial \psi}{\partial
  t}\Bigr)\,dx\,dt.
\end{align*}
Since $\mathbf{\Delta_h}u =\mathbf{\Delta_h}\psi$ when $t=0$,  the
above  expression is majorized by
\begin{gather*}
\frac{1}{2}\int_{\Omega}\zeta^p ((\mathbf{\Delta_h}\psi)_{T}^2
-(\mathbf{\Delta_h}\psi)_{0}^2)\,dx
 +\frac{1}{2}\int_0^T\!\! \int_{\Omega}\zeta^p\Bigl((\mathbf{\Delta_h}u)^2 
+\bigl(\mathbf{\Delta_h}\frac{\partial \psi}{\partial t}\bigr)^2\Bigr)\,dx\,dt,
\end{gather*}
where the inequality $2 \mathbf{\Delta_h}u\,\mathbf{\Delta_h}\psi \leq
(\mathbf{\Delta_h}u)^2 + (\mathbf{\Delta_h}\psi)^2$ was used at time
$T$.

At this stage there are no ``forbidden'' time derivatives left and so
we may safely let $\varepsilon$ go to zero. By Lemma 3 we may pass to
the limit under the integral sign and hence  the estimate for the limit $u$
(no longer $u^{\varepsilon}$) becomes
\begin{align*}
&\int_0^T\!\!\int_{\Omega}(\langle\mathbf{\Delta_h} \mathbf{A},\nabla
(\zeta^p\mathbf{\Delta_h}u)\rangle \,dx\,dt\\
\leq& \int_0^T\!\!\int_{\Omega}(\langle\mathbf{\Delta_h} \mathbf{A},\nabla
(\zeta^p\mathbf{\Delta_h}\psi)\rangle \,dx\,dt\\
+&\frac{1}{2}\int_0^T\!\! \int_{\Omega}\zeta^p \Bigl((\mathbf{\Delta_h}u)^2 
+\bigl(\mathbf{\Delta_h}\frac{\partial \psi}{\partial t}\bigr)^2\Bigr)\,dx\,dt +
 \frac{1}{2}\int_{\Omega}\zeta^p (\mathbf{\Delta_h}\psi)_{T}^2\,dx,
\end{align*}
where $$\mathbf{\Delta_h} \mathbf{A} = \mathbf{A}(x+h,t) - \mathbf{A}(x,t).$$
We write this more conveniently as
\begin{align}
&\int_0^T\!\!\int_{\Omega}\zeta^p(\langle\mathbf{\Delta_h} \mathbf{A},\nabla
\mathbf{\Delta_h}u)\rangle \,dx\,dt \nonumber\\
\leq& \overbrace{\int_0^T\!\!\int_{\Omega}p\zeta^{p-1}|\mathbf{\Delta_h}
\mathbf{A}||\mathbf{\Delta_h}u||\nabla \zeta|\,dx\,dt}^{\mathrm{I}} \nonumber\\
+& \overbrace{\int_0^T\!\!\int_{\Omega}p\zeta^{p-1}|\mathbf{\Delta_h}
\mathbf{A}||\mathbf{\Delta_h}\psi||\nabla \zeta|\,dx\,dt}^{\mathrm{II}}\\
+&\overbrace{\int_0^T\!\!\int_{\Omega}\zeta^{p}|\mathbf{\Delta_h}
\mathbf{A}||\nabla \mathbf{\Delta_h}\psi|\,dx\,dt}^{\mathrm{III}}\nonumber\\
+&\frac{1}{2}\int_0^T\!\! \int_{\Omega}\zeta^p\Bigl((\mathbf{\Delta_h}u)^2 
\bigl(\mathbf{\Delta_h}\frac{\partial \psi}{\partial t}\bigr)^2\Bigr)\,dx\,dt +
 \frac{1}{2}\int_{\Omega}\zeta^p (\mathbf{\Delta_h}\psi)_{T}^2\,dx.\nonumber
\end{align}
The integrand on left-hand side is $\langle\mathbf{\Delta_h} \mathbf{A},\nabla
\mathbf{\Delta_h}u)\rangle =$
\begin{gather}
\langle|\nabla u(x\!+\!h,t)|^{\frac{p-2}{2}}\nabla u(x\!+\!h,t) -|\nabla
u(x,t)|^{\frac{p-2}{2}}\nabla u(x,t), \nabla u(x\!+\!h,t)-\nabla
u(x,t)\rangle \nonumber \\
\geq \frac{4}{p^2}\left|\mathbf{F}(x+h,t) -\mathbf{F}(x,t)\right|^2
= \frac{4}{p^2}|\mathbf{\Delta_h}\mathbf{F}|^{2},
\end{gather}
where the elementary inequality 
$$\frac{4}{p^2}\left||b|^{\frac{p-2}{2}}b-|a|^{\frac{p-2}{2}}a\right|^2
\leq \langle|b|^{p-2} b-|a|^{p-2}a,b-a\rangle$$
for vectors was used. We aim at  an estimate for the integral
of $\zeta^p |\mathbf{\Delta_h}\mathbf{F}|²$.
 
We divide the $\mathbf{\Delta_h}$-terms  by $|h|$ so that the desired
difference quotients appear. The estimate
$$\left|\frac{\mathbf{\Delta_h}\mathbf{A}}{h}\right| \leq (p-1)\left| 
\frac{\mathbf{\Delta_h}\mathbf{F}}{h}\right|\bigl(|\nabla
u(x+h,t)|^{\frac{p-2}{2}}+|\nabla u(x,t)|^{\frac{p-2}{2}}\bigr),$$
coming from the elementary vector inequality
\begin{gather*}
\left||b|^{p-2} b-|a|^{p-2}a\right|
\leq(p-1)\Bigl(|b|^{\frac{p-2}{2}}+|a|^{\frac{p-2}{2}}\Bigr)\left||b|^{\frac{p-2}{2}}b-|a|^{\frac{p-2}{2}}a\right|,
\end{gather*}
is used in the integrands of I, II, and III. In I we split the factors so
that
\begin{gather*}
p\zeta^{p-1}\left|\frac{\mathbf{\Delta_h}\mathbf{A}}{h}\right|\left|\frac{\mathbf{\Delta_h}u}{h}\right||\nabla
\zeta|\\
\leq
p(p-1)\!\left[\zeta^{\frac{p}{2}}\bigl\vert
  \frac{\mathbf{\Delta_h}\mathbf{F}}{h}
\bigr\vert \right]\left[\bigl\vert \frac{\mathbf{\Delta_h}u}{h}\bigr\vert|\nabla
  \zeta|\right]\left[\zeta^{\frac{p-2}{2}}\bigl(|\nabla
    u(x,t)|^{\frac{p-2}{2}} +|\nabla
    u(x\!+\!h,t)|^{\frac{p-2}{2}}\bigr)\right]
\end{gather*}
and use Young's inequality 
$$abc \leq \frac{\varepsilon^2 a^2}{2} + \frac{\varepsilon^{-p} b^p}{p}
+  \frac{(p-2)c^{\frac{2p}{p-2}}}{2p}$$
to get the bound 
\begin{gather*}
\frac{\mathrm{I}}{|h|^{2}} \leq
\frac{p(p-1)\varepsilon^{2}}{2}\int_0^T\!\!\int_{\Omega}\zeta^p
\left|\frac{\mathbf{\Delta_h}\mathbf{F}}{h}\right|^2\,dx\,dt\\
+(p-1)\varepsilon^{-p}\int_0^T\!\!\int_{\Omega}\left|\frac{\mathbf{\Delta_h}u}{h}\right|^p|\nabla
\zeta|^p\,dx\,dt\\
+c_p \int_0^T\!\!\int_{\Omega}\zeta^p\Big(|\nabla u(x,t)|^p + |\nabla
u(x+h,t)|^p\Big) \,dx\,dt
\end{gather*}
The integral $\mathrm{II}/|h|^2$ has a similar majorant, the only difference
being that $\mathbf{\Delta_h}u$ be replaced by
$\mathbf{\Delta_h}\psi.$ The integrand of III is estimated in a
similar way:
 \begin{gather*}
p\zeta^{p}\left|\frac{\mathbf{\Delta_h}\mathbf{A}}{h}\right|\left|
\frac{\mathbf{\Delta_h}\psi}{h}\right|\\
\leq
p(p-1)\!\left[\zeta^{\frac{p}{2}}\bigl\vert
  \frac{\mathbf{\Delta_h}\mathbf{F}}{h}\bigr\vert
\right]\left[\zeta\bigl\vert
  \nabla(\frac{\mathbf{\Delta_h}u}{h})\bigr\vert \right]\left[\zeta^{\frac{p-2}{2}}\bigl(|\nabla
    u(x,t)|^{\frac{p-2}{2}} +|\nabla
    u(x\!+\!h,t)|^{\frac{p-2}{2}}\bigr)\right]\\
\leq
\frac{p(p-1)\varepsilon^2}{2}\zeta^p\left|
\frac{\mathbf{\Delta_h}\mathbf{F}}{h}\right|^2
+ (p-1)\varepsilon^{-p}\zeta^p\left|\nabla \!
  \left(\frac{\mathbf{\Delta_h}\psi}{h}\right)\right|^p\\
+c_p\zeta^p\bigl(|\nabla
    u(x,t)|^p +|\nabla
    u(x+h,t)|^p\bigr).
\end{gather*}
Adding up the three integrated estimates, we arrive at
\begin{align*}
&\qquad \qquad \frac{\mathrm{I+II+III}}{|h|^{2}}\quad  \leq\\
&\,3\,\frac{p(p-1)\varepsilon^{2}}{2}\int_0^T\!\!\int_{\Omega}\zeta^p
\left|\frac{\mathbf{\Delta_h}\mathbf{F}}{h}\right|^2\,dx\,dt\\
&+(p-1)\varepsilon^{-p}\int_0^T\!\!\int_{\Omega}\Big(\left|\frac{\mathbf{\Delta_h}u}{h}\right|^p|\nabla
\zeta|^p+ \left|\frac{\mathbf{\Delta_h}\psi}{h}\right|^p|\nabla \zeta|^p
+ \zeta^p\left|\frac{\nabla
    \big(\mathbf{\Delta_h}\psi}{h}\right|^p \Big)\,dx\,dt\\
&+3 c_p \int_0^T\!\!\int_{\Omega}\zeta^p\Big(|\nabla u(x,t)|^p + |\nabla
u(x+h,t)|^p\Big) \,dx\,dt.
\end{align*}

This complements (9). Recall (10). The next step is to
absorb
 the first integral above in the right-hand member
into the minorant in (10) by fixing $\varepsilon$ small enough, say
$$3\,\frac{p(p-1)\,\varepsilon^{2}}{2} = \frac{2}{p^2}.$$ The resulting
estimate, written out without abbreviations, is
\begin{align*}
\int_0^T\!\!\int_{\Omega}&\zeta^p\left|\frac{\mathbf{F}(x+h,t)-\mathbf{F}(x,t)}{h}\right|^{2}\,dx\,dt\\
 \leq &a_p\int_0^T\!\!\int_{\Omega}\left|\frac{u(x+h,t)-u(x,t)}{h}\right|^{p}|\nabla
\zeta|^{p}\,dx\,dt\\
+&a_p\int_0^T\!\!\int_{\Omega}\left|\frac{\psi(x+h,t)-\psi(x,t)}{h}\right|^{p}|\nabla
\zeta|^{p}\,dx\,dt\\
+&a_p\int_0^T\!\!\int_{\Omega}\zeta^{p}\left|\frac{\nabla
    \psi(x+h,t)-\nabla \psi(x,t)}{h}\right|^{p}\,dx\,dt\\
+&b_p\int_0^T\!\!\int_{\Omega}\zeta^p\Big(|\nabla u(x,t)|^p + |\nabla
u(x+h,t)|^p\Big) \,dx\,dt\\
+&c_p\int_0^T\!\!\int_{\Omega}\zeta^{p}\left|\frac{u(x+h,t)-u(x,t)}{h}\right|^{2}\,dx\,dt\\
+&c_p\int_0^T\!\!\int_{\Omega}\zeta^{p}\left|\frac{\psi_t(x+h,t)-\psi_t(x,t)}{h}\right|^{2}\,dx\,dt\\
+&c_p
\int_{\Omega}\zeta^{p}\left|\frac{\psi(x+h,T)-\psi(x,T)}{h}\right|^{2}\,dx,
\end{align*}
where the constants depend only on $p$. Finally, letting the increment
$h \rightarrow 0$ in any desired direction, we arrive at the
estimate in the theorem. Here we use the characterization of Sobolev
spaces in terms of integrated differential quotients,
cf. [G, Chapter 8.1]. This concludes our proof of Theorem 5.

\begin{cor} If $u$ is the solution to the obstacle problem with the
  obstacle $\psi$, then $\Delta_pu = \nabla \cdot (|\nabla
  u|^{p-2}\nabla u)$ belongs to the space
  $L^{\frac{p}{p-1}}_{loc}(\Omega_T)$ and
$$\int_0^T\!\!\int_{\Omega}\varphi\,\Delta_p u \,dx\,dt =
\int_0^T\!\!\int_{\Omega}\langle \, |\nabla
  u|^{p-2}\nabla u,\nabla \varphi \rangle \,dx\,dt$$
for all test functions $\varphi$ in $C_0^{\infty}(\Omega_T).$
\end{cor}

\emph{Proof:} Since $\mathbf{F}$ is in Sobolev's space and p > 2, we can
differentiate
$$|\nabla u|^{p-2}\nabla u = |\mathbf{F}|^{\frac{p-2}{p}}\mathbf{F}$$
and hence
$$\left|\frac{\partial}{\partial x_j}(|\nabla u|^{p-2}\nabla u)\right|
\leq 2\Bigl(1-\frac{1}{p}\Bigr)
 |\mathbf{F}|^{\frac{p-2}{p}}\left|\frac{\partial
    \mathbf{F}}{\partial x_j} \right|.$$
By H\"{o}lder's inequality 
$$ \frac{\partial}{\partial x_j}(|\nabla u|^{p-2}\nabla u) \in
L^{\frac{p}{p-1}}_{loc}(\Omega_T),$$
since $\mathbf{F} \in L^{2}(\Omega_T)$ and $D\mathbf{F} \in
L^{2}(\Omega_T).$ $\Box$

\section{The Time Derivative}

For the proof of the Theorem we notice that the contact set $\Xi =
\{u = \psi\}$ is  a closed subset of $\overline{\Omega_T}$ and that
 its complement $\Upsilon =\Omega_T \setminus
\Xi$ is open. In the set $\Upsilon$, where the obstacle does not hinder, $u$
is a solution to the Evolutionary $p$-Laplace Equation $u_t =
\Delta_p u.$
In other words, whenever $\phi \in C_0^{\infty}(\Upsilon)$,
$$ \int_0^T\!\!\int_{\Omega} u \frac{\partial \phi}{\partial
  t}\,dx\,dt = \int_0^T\!\!\int_{\Omega}\langle \, |\nabla
  u|^{p-2}\nabla u,\nabla \phi \rangle \,dx\,dt = -
  \int_0^T\!\!\int_{\Omega}\phi\Delta_p u \,dx\,dt,$$
the actual set of integration being $\Upsilon.$ Here Corollary was used. Thus $u_t$ is
 available, but only in
$\Upsilon$ to begin with. (See also [L2].). Let $\phi$ denote
an arbitrary test function in $C_0^{\infty}(\Omega_T).$ We need a specific test
function with compact support in $\Upsilon$. To construct it, define
$$\theta_k = \min\{1,k(u-\psi)\}, \qquad k = 1,2,\cdots.$$
Then $1-\theta_k = 1$ in $\Xi$ and pointwise the monotone convergence
$1-\theta_k
\longrightarrow \chi_{\Xi}$ holds. Moreover, the support of $\theta_k$ is
compact in $\Upsilon$. The time derivative of $\theta_k$ is available! 
 
Using
$$ \int_0^T\!\!\int_{\Omega} u \frac{\partial }{\partial
  t}(\theta_k \phi)\,dx\,dt = \int_0^T\!\!\int_{\Omega}\langle \, |\nabla
  u|^{p-2}\nabla u,\nabla(\theta_k \phi) \rangle \,dx\,dt, $$
we write
\begin{align*}
&\int_0^T\!\!\int_{\Omega}\phi\Delta_p u \,dx\,dt =
  -\int_0^T\!\!\int_{\Omega}\langle \, |\nabla
  u|^{p-2}\nabla u,\nabla \phi \rangle \,dx\,dt\\
= -&\int_0^T\!\!\int_{\Omega}\langle \, |\nabla
  u|^{p-2}\nabla u,\nabla (\theta_k \phi + (1-\theta_k) \phi)  \rangle
  \,dx\,dt\\
= -&\int_0^T\!\!\int_{\Omega}\langle \, |\nabla
  u|^{p-2}\nabla u,\nabla (\theta_k \phi)  \rangle
  \,dx\,dt -\int_0^T\!\!\int_{\Omega}\langle \, |\nabla
  u|^{p-2}\nabla u,\nabla \left((1-\theta_k \right) \phi)  \rangle\,dx\,dt \\
=  -&\int_0^T\!\!\int_{\Omega} u \frac{\partial }{\partial
  t}(\theta_k \phi)\,dx\,dt
 + \int_0^T\!\!\int_{\Omega}(1-\theta_k) \phi \Delta_p u
  \,dx\,dt. \\
\end{align*}
The last integral has the limit
$$ \lim_{k \rightarrow 0}\, \int_0^T\!\!\int_{\Omega}(1-\theta_k) \phi\, \Delta_p u
  \,dx\,dt = \underset{\Xi\phantom{ab}}{\int\!\!\int}\phi\,\Delta_p \psi
  \,dx\,dt.$$
In the integral with the time derivative we write
$$-u\frac{\partial }{\partial t}(\theta_k \phi) =
 -u\frac{\partial \phi }{\partial t} +(u-\psi)\frac{\partial }{\partial
  t}( 1-\theta_k) \phi + \psi \frac{\partial }{\partial
  t}( 1-\theta_k) \phi $$
and obtain
\begin{align*}
 -&\int_0^T\!\!\int_{\Omega} u \frac{\partial }{\partial
  t}(\theta_k \phi)\,dx\,dt = \int_0^T\!\!\int_{\Omega}
 -u\frac{\partial \phi }{\partial t}\,dx\,dt\\
 +
 &\int_0^T\!\!\int_{\Omega}
(u-\psi)\frac{\partial }{\partial
  t}(( 1-\theta_k) \phi) \,dx\,dt - \int_0^T\!\!\int_{\Omega}(
1-\theta_k) \phi \frac{\partial \psi }{\partial t}\,dx\,dt,
\end{align*}
where an integration by parts has produced the last integral. It has
the evident limit
$$\lim_{k \rightarrow 0}\, \int_0^T\!\!\int_{\Omega}(1-\theta_k)
 \phi \frac{\partial \psi }{\partial t}\,dx\,dt =
\underset{\Xi\phantom{ab}}{\int\!\!\int}\phi\,\frac{\partial \psi }{\partial t}
  \,dx\,dt.$$
The middle integral vanishes as $k \rightarrow 0$:
\begin{align*}
  &\int_0^T\!\!\int_{\Omega}(u-\psi)\frac{\partial }{\partial
  t}(( 1-\theta_k) \phi) \,dx\,dt\\
 =& \int_0^T\!\!\int_{\Omega}(u-\psi)( 1-\theta_k)
 \frac{\partial \phi }{\partial t} \,dx\,dt -
 \int_0^T\!\!\int_{\Omega} \phi(u-\psi)\frac{\partial \theta_k}{\partial
  t} \,dx\,dt \\
 =& \int_0^T\!\!\int_{\Omega}(u-\psi)( 1-\theta_k)
 \frac{\partial \phi }{\partial t} \,dx\,dt -
 \frac{1}{2k} \int_0^T\!\!\int_{\Omega} \phi \frac{\partial }{\partial
  t}\theta_{k}^{2} \,dx\,dt \\
 =& \int_0^T\!\!\int_{\Omega}(u-\psi)( 1-\theta_k)
 \frac{\partial \phi }{\partial t} \,dx\,dt +
 \frac{1}{2k} \int_0^T\!\!\int_{\Omega} \theta_{k}^{2}\,
  \frac{\partial \phi }{\partial
  t}\,dx\,dt \quad \underset{k \rightarrow \infty}{\longrightarrow} \quad 0\,+\,0.
\end{align*}

Collecting results,
\begin{gather*}
\int_0^T\!\!\int_{\Omega}\phi \,\Delta_pu\,dx\,dt =
 -\int_0^T\!\!\int_{\Omega}\phi_t\,dx\,dt -
 \underset{\Xi\phantom{ab}}{\int\!\!\int}(\psi_t -
 \Delta_p \psi)\phi\,dx\,dt.
\end{gather*}
In other words, the final formula
\begin{gather*}
 - \int_0^T\!\!\int_{\Omega}u\phi_t \,dx\,dt = 
\int_0^T\!\!\int_{\Omega}\phi\,[\Delta_pu +
 (\psi_t- \Delta_p\psi)\chi_{\Xi}]\,dx\,dt
\end{gather*}
holds for every $\phi$ in $C_{0}^{\infty}(\Omega_T)$. Therefore
$$u_t = \Delta_pu +
 (\psi_t- \Delta_p\psi)\chi_{\Xi}$$
and this is a \emph{function} belonging to
$L^{p/(p-1)}_{loc}(\Omega_T)$.
This concludes the proof of Theorem 2. $\Box$

\bigskip
{\small {\textsf{Peter Lindqvist\\ Department of
    Mathematical Sciences, Norwegian University of Science and
    Techno\-logy, N--7491, Trondheim, Norway}}

\end{document}